\documentclass[10pt,english]{amsart}
\usepackage{helvet}
\usepackage[T1]{fontenc}
\usepackage[latin1]{inputenc}
\usepackage{a4wide}
\usepackage{amssymb}

 \theoremstyle{plain}    
 \newtheorem{thm}{Theorem}[section]
 \numberwithin{equation}{section} 
 \numberwithin{figure}{section} 
 \theoremstyle{plain}    
 \newtheorem*{thm*}{Theorem} 
 \theoremstyle{definition}
  \newtheorem*{example*}{Example}
 \theoremstyle{remark}    
 \newtheorem{notation}[thm]{Notation} 
 \theoremstyle{definition}
 \newtheorem{defn}[thm]{Definition}
 \theoremstyle{definition}
  \newtheorem{example}[thm]{Example}
 \theoremstyle{remark}
 \newtheorem{rem}[thm]{Remark}
 \theoremstyle{plain}    
 \newtheorem{prop}[thm]{Proposition} 
 \theoremstyle{plain}    
 \newtheorem{cor}[thm]{Corollary} 
 \theoremstyle{plain}    
 \newtheorem{lem}[thm]{Lemma} 

\usepackage{latexsym,amsfonts,amssymb}

\theoremstyle{remark}

\theoremstyle{definition}
\newtheorem{enavant}[thm]{}

\usepackage{babel}

\begin{document}

\title{Additive group actions on Danielewski varieties and the cancellation
problem}

\maketitle
 \insert\footins{\footnotesize\noindent \textbf{Mathematics Subject Classification (2000)}: 14R10,14R20.

\par}

\begin{center}\author{{\bf A. DUBOULOZ} \\ \vspace{0.5cm}   \address Institut Fourier, Laboratoire de Mathématiques\\ UMR5582 (UJF-CNRS)\\ BP 74, 38402 ST MARTIN D'HERES CEDEX\\ FRANCE \\ \email adrien.dubouloz@ujf-grenoble.fr \\ \vspace{0.8cm}} \end{center}

\begin{abstract}
The cancellation problem asks if two complex algebraic varieties $X$
and $Y$ of the same dimension such that $X\times\mathbb{C}$ and
$Y\times\mathbb{C}$ are isomorphic are isomorphic. Iitaka and Fujita
\cite{IiF77} established that the answer is positive for a large
class of varieties of any dimension. In 1989, Danielewski \cite{Dan89}
constructed a famous counter-example using smooth affine surfaces
with additive group actions. His construction was further generalized
by Fieseler \cite{Fie94} and Wilkens \cite{Wil98} to describe a
larger class of affine surfaces. Here we construct higher dimensional
analogues of these surfaces. We study algebraic actions of the additive
group $\mathbb{C}_{+}$ on certain of these varieties, and we obtain
counter-examples to the cancellation problem in any dimension $n\geq2$. 
\end{abstract}

\section*{Introduction}

The Cancellation Problem, which is sometimes referred to a Zariski's
Problem although Zariski's original question was different, has been
already discussed in the early seventies as the question of uniqueness
of coefficients rings. The problem at that time was to decide for
which rings $A$ and $B$ an isomorphism of the polynomials rings
$A\left[x\right]$ and $B\left[x\right]$ implies that $A$ and $B$
are isomorphic (see e.g. \cite{EaH73}). Using the fact that the tangent
bundle of the real $n$-sphere is stably trivial but not trivial,
Hochster \cite{Ho72} showed that this fails in general. 

A geometric formulation of the Cancellation Problem asks if two algebraic
varieties $X$ and $Y$ such that $Y\times\mathbb{A}^{1}$ is isomorphic
to $X\times\mathbb{A}^{1}$ are isomorphic. Clearly, if either $X$
or $Y$ does not contain rational curves, for instance $X$ or $Y$
is an abelian variety, then every isomorphism $\Phi:X\times\mathbb{A}\stackrel{\sim}{\rightarrow}Y\times\mathbb{A}^{1}$
induces an isomorphism of $X$ and $Y$. So the Cancellation Problem
leads to decide if a given algebraic variety $X$ contains a family
of rational curves, where by a rational curve we mean the image of
a nonconstant morphism $f:C\rightarrow X$, where $C$ is isomorphic
to $\mathbb{A}^{1}$ or $\mathbb{P}^{1}$. Iitaka and Fujita \cite{IiF77}
carried a geometric attack to this question using ideas from the classification
theory of noncomplete varieties. Every complex algebraic variety $X$
embeds as an open subset of complete variety $\bar{X}$ for which
the boundary $D=\bar{X}\setminus X$ is a divisor with normal crossing.
By replacing the usual sheaves of forms $\Omega^{q}\left(\bar{X}\right)$
on $\bar{X}$ by the sheaves $\Omega^{q}\left(\log D\right)$ of rational
$q$-forms having at worse logarithmic poles along $D$, Iitaka \cite{Ii77}
introduced, among others invariants, the notion of logarithmic Kodaira
dimension $\bar{\kappa}\left(X\right)$ of a noncomplete variety $X$,
which is an analogue of the usual notion of Kodaira for complete varieties.
They established the following result.

\begin{thm*}
Let $X$ and $Y$ be two nonsingular algebraic varieties and assume
that either $\bar{\kappa}\left(X\right)\geq0$ or $\bar{\kappa}\left(Y\right)\geq0$.
Then every isomorphism $\Phi:X\times\mathbb{C}\stackrel{\sim}{\rightarrow}Y\times\mathbb{C}$
induces an isomorphism between $X$ and $Y$. 
\end{thm*}
\noindent The hypothesis $\bar{\kappa}\left(X\right)\geq0$ above
guarantees that $X$ cannot contain too many rational curves. For
instance, there is no cylinder-like open subset $U\simeq C\times\mathbb{A}^{1}$
in $X$, for otherwise we would have $\bar{\kappa}\left(X\right)=-\infty$%
\footnote{Actually, a nonsingular affine surface has logarithmic Kodaira dimension
$-\infty$ if and only if its contains a cylinder-like open set (see
e.g. \cite{MiyBook}). %
}. It turns out that this additional assumption is essential, as shown
by the following example due to Danielewski \cite{Dan89}.

\begin{example*}
\emph{The surfaces} $S_{1},S_{2}\subset\mathbb{C}^{3}$ \emph{with
equations} $xz-y^{2}+1=0$ \emph{and} $x^{2}z-y^{2}+1=0$ \emph{are
not isomorphic but} $S_{1}\times\mathbb{C}$ and $S_{2}\times\mathbb{C}$
\emph{are}. In the construction of Danielewski, these surfaces appear
as the total spaces of principal homogeneous $\mathbb{C}_{+}$-bundles
over $\tilde{\mathbb{A}}$, the affine line with a double origin,
obtained by identifying two copies of $\mathbb{A}^{1}$ along $\mathbb{A}^{1}\setminus\left\{ 0\right\} $.
The isomorphism $S_{1}\times\mathbb{C}\simeq S_{2}\times\mathbb{C}$
is obtained by forming the fiber product $S_{1}\times_{\tilde{\mathbb{A}}}S_{2}$,
which is a principal $\mathbb{C}_{+}$-bundle over both $S_{1}$ and
$S_{2}$, and using the fact that every such bundle over an affine
variety is trivial. On the other hand, $S_{1}$ and $S_{2}$ are not
homeomorphic when equipped with the complex topology. More precisely,
Danielewski established that the fundamental groups at infinity of
$S_{1}$ and $S_{2}$ are isomorphic to $\mathbb{Z}/2\mathbb{Z}$
and $\mathbb{Z}/4\mathbb{Z}$ respectively. Fieseler \cite{Fie94}
studied and classified algebraic $\mathbb{C}_{+}$-actions normal
affine surfaces. As a consequence of his classification, he obtained
many new examples of the same kind (see also \cite{Wil98}).
\end{example*}
Here we construct higher dimensional analogues of Danielewski's counter-example.
The paper is organized as follows. In the first section, we introduce
a natural generalization of the surfaces $S_{1}$ and $S_{2}$ above
in the form of certain affine varieties which are the total spaces
of certain principal homogeneous $\mathbb{C}_{+}$-bundle over $\tilde{\mathbb{A}}^{n}$,
the affine $n$-space with a multiple system of coordinate hyperplanes.
We call them \emph{Danielewski varieties}. For instance, for every
multi-index $\left[m\right]=\left(m_{1},\ldots,m_{n}\right)\in\mathbb{Z}_{>0}^{n}$
the nonsingular hypersurface $X_{\left[m\right]}\subset\mathbb{C}^{n+2}$
with equation $x_{1}^{m_{1}}\cdots x_{n}^{m_{n}}z=y^{2}-1$ is a Danielewski
variety. As a generalization of a result of Danielewski (see also
\cite{Fie94}), we establish that the total space of a principal homogeneous
$\mathbb{C}_{+}$-bundle over $\tilde{\mathbb{A}^{n}}$ is a Danielewski
variety if and only if it is separated. This leads to simple description
of these varieties in terms of $\check{C}$ech cocycles (see Theorem
(\ref{thm:Main-Theorem_1})). 

In a second part, we study algebraic $\mathbb{C}_{+}$-actions on
a certain class of varieties which contains the Danielewski varieties
$X_{\left[m\right]}$ as above. In particular we compute the Makar-Limanov
invariant of these varieties, \emph{i.e.} the set of regular functions
invariant under \emph{all} $\mathbb{C}_{+}$-actions. We obtain the
following generalisation of a result due to Makar-Limanov \cite{ML01}
for the case of surfaces (see Theorem (\ref{thm:Main-Theorem_2})
). 

\begin{thm*}
If $\left(m_{1},\ldots,m_{n}\right)\in\mathbb{Z}_{>1}^{n}$ then the
Makar-Limanov of a variety $X\subset\mathbb{C}^{n+2}$ with equation
\[
x_{1}^{m_{1}}\cdots x^{m_{n}}z=y^{r}+\sum a_{i}\left(x_{1},\ldots,x_{n}\right)y^{i},\quad\textrm{where }r\geq,\]
 is isomorphic to $\mathbb{C}\left[x_{1},\ldots,x_{n}\right]$. 
\end{thm*}
\noindent As a consequence, we obtain infinite families of counter-examples
to the Cancellation Problem in every dimension $n\geq2$.

\begin{thm*}
Let $\left[m\right]=\left(m_{1},\ldots,m_{n}\right)\in\mathbb{Z}_{>1}^{n}$
and $\left[m'\right]=\left(m'_{1},\ldots,m'_{n}\right)\in\mathbb{Z}_{>1}^{n}$
be two multi-index for which the subsets $\left\{ m_{1},\ldots,m_{n}\right\} $
and $\left\{ m'_{1},\ldots,m'_{n}\right\} $ of $\mathbb{Z}$ are
distint, and let $\lambda_{1},\ldots,\lambda_{r}$, where $r\geq2$
be a collection of pairwise distinct complex numbers. Then the Danielewski
varieties $X$ and $X'$ in $\mathbb{C}^{n+2}$ with equations \[
x_{1}^{m_{1}}\cdots x_{n}^{m_{n}}z-\prod_{i=1}^{r}\left(y-\lambda_{i}\right)=0\quad\textrm{and}\quad x_{1}^{m'_{1}}\cdots x_{n}^{m'_{n}}z-\prod_{i=1}^{r}\left(y-\lambda_{i}\right)=0\]
 are not isomorphic, but the varieties $X\times\mathbb{C}$ and $X'\times\mathbb{C}$
are isomorphic. 
\end{thm*}

\section{Danielewski varieties}

Danielewski's construction can be easily generalized to produce examples
of affine varieties $X$ and $Y$ such that $X\times\mathbb{C}$ and
$Y\times\mathbb{C}$ are isomorphic. Indeed, if we can equip two affine
varieties $X$ and $Y$ with structures of principal homogeneous $\mathbb{C}_{+}$-bundle
$\rho_{X}:X\rightarrow Z$ and $\rho_{Y}:Y\rightarrow Z$ over a certain
scheme $Z$, then the fiber product $X\times_{Z}Y$ will be a principal
homogeneous $\mathbb{C}_{+}$-bundle over $X$ and $Y$, whence a
trivial principal bundle $X\times\mathbb{C}\simeq X\times_{Z}Y\simeq Y\times\mathbb{C}$
as $X$ and $Y$ are both affine. The base scheme $Z$ which arises
in Danielewski's counter-example is the affine with a double origin.
The most natural generalization is to consider an \emph{affine space}
$\mathbb{C}^{n}$ \emph{with a multiple system of coordinate hyperplanes}
as a base scheme. 

\begin{notation}
In the sequel we denote the polynomial ring $\mathbb{C}\left[x_{1},\ldots,x_{n}\right]$
by $\mathbb{C}\left[\underline{x}\right]$, and the algebra $\mathbb{C}\left[x_{1},x_{1}^{-1}\ldots,x_{n},x_{n}^{-1}\right]$
of Laurent polynomials in the variables $x_{1},\ldots,x_{n}$ by $\mathbb{C}\left[\underline{x},\underline{x}^{-1}\right]$.
For every multi-index $\left[r\right]=\left(r_{1},\ldots,r_{n}\right)\in\mathbb{Z}^{n}$,
we let $\underline{x}^{\left[r\right]}=x_{1}^{r_{1}}\cdots x_{n}^{r_{n}}\in\mathbb{C}\left[\underline{x},\underline{x}^{-1}\right]$.
We denote by $H_{\underline{x}}=V\left(x_{1}\cdots x_{n}\right)$
the closed subvariety of $\mathbb{C}^{n}$ consisting of the disjoint
union of the $n$ coordinate hyperplanes. Its open complement in $\mathbb{C}^{n}$,
which is isomorphic to $\left(\mathbb{C}^{*}\right)^{n}$, will be
denoted by $U_{\underline{x}}$. 
\end{notation}
\begin{defn}
\label{def:Multiple_hyperplanes} We let $Z_{n,r}$ be the scheme
obtained by gluing $r$ copies $\delta_{i}:Z_{i}\stackrel{\sim}{\longrightarrow}\mathbb{C}^{n}$
of the affine space $\mathbb{C}^{n}=\textrm{Spec}\left(\mathbb{C}\left[x_{1},\ldots,x_{n}\right]\right)$
by the identity along $\left(\mathbb{C}^{*}\right)^{n}$. We call
$Z_{n,r}$ \emph{the affine} $n$\emph{-space with an $r$-fold system
of coordinate hyperplanes.} We consider it as a scheme over $\mathbb{C}^{n}$
via the morphism $\delta:Z_{n,r}\rightarrow\mathbb{C}^{n}$ restricting
to the $\delta_{i}$'s on the canonical open subset $Z_{i}$ of $Z_{n,r}$,
$i=1,\ldots,n$. 
\end{defn}
\begin{enavant} We recall that a principal homogeneous $\mathbb{C}_{+}$-bundle
over a base scheme $S$ is an $S$-scheme $\rho:X\rightarrow S$ equipped
with an algebraic action of the additive group $\mathbb{C}_{+}$,
such that there exists an open covering $\mathcal{U}=\left(S_{i}\right)_{i\in I}$
of $S$ for which $\rho^{-1}\left(S_{i}\right)$ is equivariantly
isomorphic to $S_{i}\times\mathbb{C}$, where $\mathbb{C}_{+}$ acts
by translations on the second factor, for every $i\in I$. In particular,
the total space of a principal homogeneous $\mathbb{C}_{+}$-bundle
has the structure of an $\mathbb{A}^{1}$-bundle over $S$. The set
$H^{1}\left(S,\mathbb{C}_{+}\right)$ of isomorphism classes of principal
homogeneous $\mathbb{C}_{+}$-bundles over $S$ is isomorphic to the
first cohomology group $\check{H}^{1}\left(S,\mathcal{O}_{S}\right)\simeq H^{1}\left(S,\mathcal{O}_{S}\right)$. 

\end{enavant}

\begin{defn}
\label{def:Danielewski-variety} A \emph{Danielewski variety} is a
nonsingular affine variety of dimension $n\geq2$ which is the total
space $\rho:X\rightarrow Z_{n,r}$ of a principal homogeneous $\mathbb{C}_{+}$-bundle
over $Z_{n,r}$ for a certain $r\geq1$. 
\end{defn}
\begin{example}
\label{exa:ML-surfaces} The Danielewski surfaces $S_{1}=\left\{ xz-y^{2}+1=0\right\} $
and $S_{2}=\left\{ x^{2}z-y^{2}+1=0\right\} $ above are Danielewski
varieties. Indeed, the projections $pr_{x}:S_{i}\rightarrow\mathbb{C}$,
$i=1,2$, factor through structural morphisms $\rho_{i}:S_{i}\rightarrow Z_{2,1}$
of principal $\mathbb{C}_{+}$-bundles over the affine line with a
double origin. More generally, the Makar-Limanov surfaces $S\subset\mathbb{C}^{3}$
with equations $x^{n}z-Q\left(x,y\right)=0$, where $n\geq1$ and
$Q\left(x,y\right)$ is a monic polynomial in $y$, such that $Q\left(0,y\right)$
has simple roots are Danielewski varieties. 
\end{example}
\begin{rem}
The scheme $Z_{n,r}$ over which a Danielewski variety $X$ becomes
the total space of a principal homogeneous $\mathbb{C}_{+}$-bundle
is unique up to isomorphism. Indeed, we have necessarily $n=\dim Z=\dim X-1$.
On the other hand, it follows from (\ref{txt:Gluing-construction})
below that $X$ is obtained by gluing $r$ copies of $\mathbb{C}^{n}\times\mathbb{C}$
along $\left(\mathbb{C}^{*}\right)^{n}\times\mathbb{C}$. So we deduce
by induction that $H_{n+1}\left(X,\mathbb{Z}\right)$ is isomorphic
to the direct sum of $r$ copies of $H_{n}\left(\left(\mathbb{C}^{*}\right)^{n}\times\mathbb{C},\mathbb{Z}\right)\simeq H_{n}\left(\left(\mathbb{C}^{*}\right)^{n},\mathbb{Z}\right)\simeq\mathbb{Z}$,
whence to $\mathbb{Z}^{r}$. Therefore, if $X$ admits another structure
of principal homogeneous $\mathbb{C}_{+}$-bundle $\rho':X\rightarrow Z_{n',r'}$
then $\left(n',r'\right)=\left(n,r\right)$. However, we want to insist
on the fact that \emph{this does not imply that the structural morphism}
$\rho:X\rightarrow Z_{n,r}$ \emph{on a Danielewski variety is unique,
even up to automorphisms of the base}. This question will be discussed
in (\ref{txt:Iso-problem}) below. 
\end{rem}
\begin{enavant} \label{txt:Gluing-construction} A principal homogeneous
$\mathbb{C}_{+}$-bundle $\rho:X\rightarrow Z_{n,r}$ becomes trivial
on the canonical open covering $\mathcal{U}$ of $Z_{n,r}$ be means
of the open subsets $Z_{i}\simeq\mathbb{C}^{n}$, $i=1,\ldots,r$
(see definition (\ref{def:Multiple_hyperplanes}) above). So there
exists a $\check{C}$ech $1$-cocycle \[
g=\left\{ g_{ij}\right\} _{i,j=1,\ldots,g}\in C^{1}\left(\mathcal{U},\mathcal{O}_{Z_{n,r}}\right)\simeq\bigoplus_{i=1}^{r}\mathbb{C}\left[\underline{x},\underline{x}^{-1}\right]\]
 representing the isomorphism class $\left[g\right]\in H^{1}\left(Z_{n,r},\mathcal{O}_{Z_{n,r}}\right)\simeq\check{H}^{1}\left(\mathcal{U},\mathcal{O}_{Z_{n,r}}\right)$
of $X$ such that $X$ is equivariantly isomorphic to the scheme obtained
by gluing $r$ copies $Z_{i}\times\mathbb{C}=\textrm{Spec}\left(\mathbb{C}\left[\underline{x}\right]\left[t_{i}\right]\right)$
of $\mathbb{C}^{n}\times\mathbb{C}$, equipped with $\mathbb{C}_{+}$-actions
by translations on the second factor, outside $H_{\underline{x}}\times\mathbb{C}\subset Z_{i}\times\mathbb{C}$
by means of the equivariant isomorphisms \[
\phi_{ij}:\left(Z_{j}\setminus H_{\underline{x}}\right)\times\mathbb{C}\stackrel{\sim}{\longrightarrow}\left(Z_{i}\setminus H_{\underline{x}}\right)\times\mathbb{C},\qquad\left(\underline{x},t_{j}\right)\mapsto\left(\underline{x},t_{j}+g_{ij}\left(\underline{x},\underline{x}^{-1}\right)\right),\quad i\neq j.\]
Since a Danielewski variety $X$ is affine, the corresponding transition
cocycle is not arbitrary. For instance, the trivial cocycle corresponds
to the trivial $\mathbb{C}_{+}$-bundle $Z_{n,r}\times\mathbb{C}$
which is not even separated if $r\geq2$. More generally, if one of
the rational functions $g_{ij}$ is regular at a point $\lambda=\left(\lambda_{1},\ldots,\lambda_{n}\right)\in H_{\underline{x}}\subset\mathbb{C}^{n}$,
then for every germ of curve $C\subset\mathbb{C}^{n}$ intersecting
$H_{\underline{x}}$ transversely in $\lambda$, $\left(\rho\circ\delta\right)^{-1}\left(C\right)\subset X$
is a nonseparated scheme. On the other hand, it follows from a very
general result of Danielewski that the total space of a principal
homogeneous $\mathbb{C}_{+}$-bundle $\rho:X\rightarrow Z_{n,2}$
defined by a cocycle $g_{12}=\underline{x}^{-\left[r\right]}a\left(\underline{x}\right)$,
where $\left[r\right]\in\mathbb{Z}_{\geq1}^{n}$, such that $\underline{x}^{\left[r\right]}\mathbb{C}\left[\underline{x}\right]+a\left(\underline{x}\right)\mathbb{C}\left[\underline{x}\right]=\mathbb{C}\left[\underline{x}\right]$
is affine, isomorphic to the variety $X\subset\mathbb{C}^{n+2}$ with
equation $\underline{x}^{\left[r\right]}z-y^{2}-a\left(\underline{x}\right)y=0$.
More generally, we have the following result. 

\end{enavant}

\begin{thm}
\label{thm:Bundle-charac} For the total space of a principal $\mathbb{C}_{+}$-bundle
$\rho:X\rightarrow Z_{n,r}$ defined by a transition cocycle $g=\left\{ g_{ij}\left(\underline{x},\underline{x}^{-1}\right)\right\} _{i,j=1,\ldots,r}$
the following are equivalent.

(1) For every $i\neq j$, $g_{ij}=x^{-\left[m_{ij}\right]}a_{ij}\left(\underline{x}\right)$
for a certain multi-index $\left[m_{ij}\right]\in\mathbb{Z}_{>0}^{n}$
and a polynomial $a_{ij}\left(\underline{x}\right)$ such that $a_{ij}\left(\underline{x}\right)\mathbb{C}\left[\underline{x}\right]+\underline{x}^{\left(1,\ldots,1\right)}\mathbb{C}\left[\underline{x}\right]=\mathbb{C}\left[\underline{x}\right]$, 

(2) $X$ is separated

(3) $X$ is affine. 
\end{thm}
\begin{proof}
We deduce from I.5.5.6 in \cite{EGAI} that $X$ is separated if and
only if $g_{ij}\in\mathbb{C}\left[\underline{x},\underline{x}^{-1}\right]$
generates $\mathbb{C}\left[\underline{x},\underline{x}^{-1}\right]$
as a $\mathbb{C}\left[\underline{x}\right]$-algebra for every $i\neq j$
. Letting $g_{ij}=\underline{x}^{-\left[m\right]}a\left(\underline{x}\right)$,
where $\left[m\right]\in\mathbb{Z}_{\geq0}^{n}$ and where $a\left(\underline{x}\right)\in\mathbb{C}\left[\underline{x}\right]$,
this is the case if and only if $\underline{x}^{-\left[m\right]}$
generates $\mathbb{C}\left[\underline{x},\underline{x}^{-1}\right]$
as a $\mathbb{C}\left[\underline{x}\right]$-algebra and $a\left(\underline{x}\right)\mathbb{C}\left[\underline{x}\right]+\underline{x}^{\left[m\right]}\mathbb{C}\left[\underline{x}\right]=\mathbb{C}\left[\underline{x}\right]$.
Indeed, the condition is sufficient as it guarantees that $\mathbb{C}\left[\underline{x},\underline{x}^{-1}\right]=\mathbb{C}\left[\underline{x}\right]\left[\underline{x}^{-\left[m\right]}\right]\subset\mathbb{C}\left[\underline{x}\right]\left[g_{ij}\right]$.
Conversely, if $\mathbb{C}\left[\underline{x},\underline{x}^{-1}\right]=\mathbb{C}\left[\underline{x}\right]\left[g_{ij}\right]$
then $g_{ij}=x^{-\left[m\right]}a\left(\underline{x}\right)$ for
a certain multi-index $\left[m\right]=\left(m_{1},\ldots,m_{n}\right)\in\mathbb{Z}_{\geq1}^{n}$
and a polynomial $a\in\mathbb{C}\left[\underline{x}\right]$ not divisible
by $x_{i}$ for every $i=1,\ldots,r$. Indeed, if there exists an
indice $i$ such that $m_{i}\leq0$ then $x_{i}^{-1}\not\in\mathbb{C}\left[\underline{x}\right]\left[g_{ij}\right]$
which contradicts our hypothesis. Furthermore, since $x^{-\left[m\right]}\in\mathbb{C}\left[\underline{x}\right]\left[g_{ij}\right]$,
there exists polynomials $b_{1},\ldots,b_{s}\in\mathbb{C}\left[\underline{x}\right]$
such that $x^{-\left[m\right]}=b_{0}+b_{1}ax^{-\left[m\right]}+\ldots+b_{s}a^{-s\left[m\right]}\in\mathbb{C}\left[\underline{x}\right]\left[g_{ij}\right]$.
This means equivalently that $\underline{x}^{\left(s-1\right)\left[m\right]}=b_{0}\underline{x}^{s\left[m\right]}+ca$
for a certain $c\in\mathbb{C}\left[\underline{x}\right]$. If $s\neq1$
then $c\in\underline{x}^{\left(s-1\right)\left[m\right]}\mathbb{C}\left[\underline{x}\right]$
as the $x_{i}$'s do not divide $a$, and so, there exists $c'\in\mathbb{C}\left[\underline{x}\right]$
such that $1=b_{0}\underline{x}^{-\left[m\right]}+c'a$. This proves
that (1) and (2) are equivalent.

Now it remains to show that if the $g_{ij}=x^{-\left[m_{ij}\right]}a_{ij}\left(\underline{x}\right)$
satisfy (1), then $X$ is affine. We first observe that there exists
an indice $i_{0}$ such that $m_{1i_{0},k}=\max\left\{ m_{1i,k}\right\} $
for every $i=2,\ldots,r$ and every $k=1,\ldots,n$. Indeed, suppose
on that there exists two indices $i\neq j$, say $i=2$ and $j=3$,
and two indices $l\neq k$ such that $m_{12,k}<m_{13,k}$ but $m_{12,l}>m_{13,l}$.
We let $\left[\mu\right]\in\mathbb{Z}_{\geq0}^{n}$ be the multi-index
with components $\mu_{s}=\max\left(m_{12,s},m_{13,s}\right)$, so
that $\mu_{k}-m_{13,k}=0$ and $\mu_{k}-m_{12,k}>0$ whereas $\mu_{l}-m_{12,l}=0$
and $\mu_{l}-m_{13,l}>0$. It follows from the cocycle relation $g_{23}=g_{13}-g_{12}$
that \begin{eqnarray*}
\underline{x}^{\left[\mu\right]-\left[m_{23}\right]}a_{23}\left(\underline{x}\right) & = & \underline{x}^{\left[\mu\right]-\left[m_{13}\right]}a_{13}\left(\underline{x}\right)-\underline{x}^{\left[\mu\right]-\left[m_{12}\right]}a_{12}\left(\underline{x}\right)\in\left(x_{k},x_{l}\right)\mathbb{C}\left[\underline{x}\right]\subset\mathbb{C}\left[\underline{x}\right].\end{eqnarray*}
 Since the $x_{i}$'s do not divide the $a_{ij}$'s, it follows that
neither $x_{k}$ nor $x_{l}$ divides the polynomial on the right.
Thus $m_{23,l}=\mu_{l}$ and $m_{23,k}=\mu_{k}$. This implies that
$a_{23}\left(\underline{x}\right)\in\left(x_{k},x_{l}\right)\mathbb{C}\left[\underline{x}\right]$
which contradicts (1) above. Therefore, the subset of $\mathbb{Z}^{n}$
consisting of the multi-indices $\left[m_{1i}\right]$, $i=2,\ldots,r$,
is totally ordered for the restriction of the product ordering of
$\mathbb{Z}^{n}$, and so, there exists an indice $i_{0}$ such that
$m_{1i_{0},k}=\max\left\{ m_{1i,k}\right\} $ for every $i=2,\ldots,r$
and every $k=1,\ldots,n$. By construction, $\sigma_{i}\left(\underline{x}\right)=\underline{x}^{\left[m_{1i_{0}}\right]}g_{1i}\left(\underline{x},\underline{x}^{-1}\right)$
is a polynomial every $i=2,\ldots,r$, and $\sigma_{i_{0}}\left(\underline{x}\right)$
restricts to a nonzero constant $\lambda\in\mathbb{C}^{*}$ on $H_{\underline{x}}\subset\mathbb{C}^{n}$.
Letting $\sigma_{1}\left(\underline{x}\right)=0$, we deduce from
the cocycle relation that $\underline{x}^{\left[m_{1i_{0}}\right]}g_{ij}=\left(\sigma_{j}\left(\underline{x}\right)-\sigma_{i}\left(\underline{x}\right)\right)$
for every $i\neq j$. In turn, this implies that the local morphisms
\[
\psi_{i}:Z_{i}\times\mathbb{C}=\textrm{Spec}\left(\mathbb{C}\left[\underline{x}\right]\left[t_{i}\right]\right)\longrightarrow\mathbb{C}^{n}\times\mathbb{C},\quad\left(\underline{x},t_{i}\right)\mapsto\left(\underline{x},\underline{x}^{\left[m_{1i_{0}}\right]}t_{i}+\sigma_{i}\left(\underline{x}\right)\right),\quad i=1,\ldots,r\]
 glue to a birational morphism $\psi:X\rightarrow\mathbb{C}^{n}\times\mathbb{C}$.
By construction, the images by $\psi$ of $H_{\underline{x}}\times\mathbb{C}\subset Z_{i_{0}}\times\mathbb{C}$
and $H_{\underline{x}}\times\mathbb{C}\subset Z_{1}\times\mathbb{C}$
are disjoint, contained respectively in the closed subsets $V\left(\underline{x},t-\lambda\right)$
and $V\left(\underline{x},t\right)$ of $\mathbb{C}^{n}\times\mathbb{C}=\textrm{Spec}\left(\mathbb{C}\left[\underline{x}\right]\left[t\right]\right)$.
Therefore, $\psi^{-1}\left(\mathbb{C}^{n}\times\mathbb{C}\setminus V\left(\underline{x},t\right)\right)$
is contained in the complement $V_{1}$ in $X$ of $H_{\underline{x}}\times\mathbb{C}\subset Z_{1}\times\mathbb{C}$,
whereas $\psi^{-1}\left(\mathbb{C}^{n}\times\mathbb{C}\setminus V\left(\underline{x},t-\lambda\right)\right)$
is contained in the complement $V_{i_{0}}$ in $X$ of $H_{\underline{x}}\times\mathbb{C}\subset Z_{i_{0}}\times\mathbb{C}$.
Clearly, $\rho:X\rightarrow Z_{n,r}$ restricts on $V_{1}$ and $V_{i_{0}}$
to the structural morphisms $\rho_{1}:V_{1}\rightarrow Z_{n,r-1}$
and $\rho_{i_{0}}:V_{i_{0}}\rightarrow Z_{n,r-1}$ of the principal
homogeneous $\mathbb{C}_{+}$-bundles corresponding tho the $\check{C}$ech
cocycles $\left\{ g_{ij}\right\} _{i,j=2,\ldots,r}$ and $\left\{ g_{ij}\right\} _{i,j\neq i_{0},i,j=1,\ldots,r}$.
So we conclude by a similar induction argument as in Proposition 1.4
in \cite{Fie94} that $V_{1}$ and $V_{i_{0}}$ are affine. In turn,
this implies that $\psi:X\rightarrow\mathbb{C}^{n}\times\mathbb{C}$
is an affine morphism, and so, $X$ is affine. 
\end{proof}
\noindent The following example introduce a class of Danielewski
varieties, which contains for instance the Makar-Limanov surfaces
of example (\ref{exa:ML-surfaces}). 

\begin{example}
\label{exa:Main-example} Suppose given a collection $\sigma$ of
polynomials $\sigma_{i}\left(\underline{x}\right)\in\mathbb{C}\left[\underline{x}\right]$,
$i=1,\ldots,r$, with the following properties. 

(1) $\sigma_{i}\left(0,\ldots,0\right)\neq\sigma_{j}\left(0,\ldots,0\right)$
for every $i\neq j$, 

(2) $\sigma_{i}\left(\underline{x}\right)-\sigma_{i}\left(0,\ldots,0\right)\in\underline{x}^{\left(1,\ldots,1\right)}\mathbb{C}\left[\underline{x}\right]$
for every $i=1,\ldots,r$.

\noindent Then for every multi-index $\left[m\right]=\left(m_{1},\ldots,m_{n}\right)\in\mathbb{Z}_{>0}^{n}$
the variety $X_{\left[m\right],\sigma}\subset\mathbb{C}^{n+2}$ with
equation \[
\underline{x}^{\left[m\right]}z-\prod_{i=1}^{r}\left(y-\sigma_{i}\left(\underline{x}\right)\right)=0\]
 is a Danielewski variety. 
\end{example}
\begin{proof}
Similarly as the Danielewski surfaces, a variety $X_{\left[m\right],\sigma}$
comes naturally equipped with a surjective morphism $\pi=pr_{\underline{x}}:X_{\left[m\right],\sigma}\rightarrow\mathbb{C}^{n}$,
$\left(\underline{x},y,z\right)\mapsto\underline{x}$ restricting
to a trivial $\mathbb{A}^{1}$-bundle $\pi^{-1}\left(\left(\mathbb{C}^{*}\right)^{n}\right)\simeq\left(\mathbb{C}^{*}\right)^{n}\times\mathbb{C}$
over $U_{\underline{x}}=\left(\mathbb{C}^{*}\right)^{n}$, with coordinate
$y$ on the second factor. On the other hand, it follows from our
assumptions that the fiber \[
\pi^{-1}\left(H_{\underline{x}}\right)\simeq\textrm{Spec}\left(\mathbb{C}\left[\underline{x},y,z\right]/\left(\underline{x}^{\left(1,\ldots,1\right)},\underline{x}^{\left[m\right]}z-\prod_{i=1}^{r}\left(y-\sigma_{i}\left(\underline{x}\right)\right)\right)\right)\]
 decomposes as the disjoint union of $r$ copies $D_{i}$ of $H_{\underline{x}}\times\mathbb{C}$,
with equations $\left\{ x_{1}\cdots x_{n}=0,y=\sigma_{i}\left(0\right)\right\} $,
and with coordinate $z$ on the second factor. The open subsets $\pi^{-1}\left(U_{\underline{x}}\right)\cup C_{i}$
of $X_{\left[m\right],\sigma}$ are isomorphic to $\mathbb{C}^{n}\times\mathbb{C}$
with natural coordinates $\underline{x}$ and \[
t_{i}={\displaystyle \frac{y-\sigma_{i}\left(\underline{x}\right)}{\underline{x}^{\left[r\right]}}=\frac{z}{{\displaystyle \prod_{j\neq i}\left(y-\sigma_{j}\left(\underline{x}\right)\right)}}},\quad i=1,\ldots,r,\]
 and so, $X_{\left[m\right],\sigma}$ is isomorphic to the total space
of the principal homogeneous $\mathbb{C}_{+}$-bundle defined by the
transition cocycles $g_{ij}=\underline{x}^{-\left[r\right]}\left(\sigma_{j}\left(\underline{x}\right)-\sigma_{i}\left(\underline{x}\right)\right)$,
$i,j=1,\ldots,r$. 
\end{proof}
\noindent As a consequence of the general principle discussed at
the beginning of this section, Danielewski varieties are natural candidates
for being counter-examples to the Cancellation problem. 

\begin{prop}
\label{pro:Isomorphic-cylinder} If two Danielewski varieties $X_{1}$
and $X_{2}$ are the total spaces of $\mathbb{C}_{+}$-principal bundles
over the same base $Z_{n,r}$ then $X_{1}\times\mathbb{C}$ and $X_{2}\times\mathbb{C}$
are isomorphic. 
\end{prop}
\begin{example}
Given a polynomial $P\left(y\right)\in\mathbb{C}\left[y\right]$ with
$r\geq2$ simple roots, the varieties $\tilde{X}_{\left[m\right],P}\subset\mathbb{C}^{n+3}=\textrm{Spec}\left(\mathbb{C}\left[\underline{x},y,z,u\right]\right)$
with equations $\underline{x}^{\left[m\right]}z-P\left(y\right)=0$,
where $\left[m\right]\in\mathbb{Z}_{\geq1}^{n}$ is an arbitrary multi-index,
are all isomorphic. Indeed $\tilde{X}_{\left[m\right],P}$ is isomorphic
to $X_{\left[m\right],P}\times\mathbb{C}$, where $X_{\left[m\right],P}\subset\mathbb{C}^{n+2}=\textrm{Spec}\left(\mathbb{C}\left[\underline{x},y,z\right]\right)$
denotes the Danielewski variety with equation $x^{\left[m\right]}z-P\left(y\right)=0$,
which has the structure of a principal homogeneous $\mathbb{C}_{+}$-bundle
over $Z_{n,r}$ (see example (\ref{exa:Main-example})). 
\end{example}
\begin{enavant} \label{txt:Iso-problem} This leads to the difficult
problem of deciding which Danielewski varieties are isomorphic as
abstract varieties. Things would be simpler if the structural morphism
$\rho:X\rightarrow Z_{n,r}$ on a Danielewski variety were unique
up to automorphisms of the base. However, this is definitely not the
case in general, as shown by the Danielewski surface $S_{1}=\left\{ xz-y^{2}+1=0\right\} \subset\mathbb{C}^{3}$,
which admits two such structures, due to the symmetry between the
variables $x$ and $z$. Actually, the situation is even worse since
in general, a Danielewski variety admitting a second $\mathbb{C}_{+}$-action,
whose general orbits are distinct from the general fibers of the structural
morphism $\rho:X\rightarrow Z_{n,r}$, comes equipped with a one parameter
family of distinct structures of principal homogeneous $\mathbb{C}_{+}$-bundles.
Indeed, let $G_{1}\simeq\mathbb{C}_{+}$ and $G_{2}\simeq\mathbb{C}_{+}$
be one-parameter subgroups of $\textrm{Aut}\left(X\right)$ corresponding
respectively to a principal homogeneous $\mathbb{C}_{+}$-bundle structure
on $\rho:X\rightarrow Z_{n,r}$ and another nontrivial $\mathbb{C}_{+}$-action
on $X$ with general orbits distinct from the ones of $G_{1}$. Then
the subgroups $\phi_{t}^{-1}G_{1}\phi_{t}\simeq\mathbb{C}_{+}$ of
$\textrm{Aut}\left(X\right)$, where $\phi_{t}\in G_{2}$, correspond
to principal homogeneous $\mathbb{C}_{+}$-bundle structures on $X$,
with pairwise distinct general orbits provided that the generators
of $G_{1}$ and $G_{2}$ do not commute. 

\end{enavant}

\begin{enavant} There exists a useful geometric criterion to decide
if a smooth affine surface admits two $\mathbb{C}_{+}$-actions with
distinct general orbits. As is well-known, there exists a correspondence
between algebraic $\mathbb{C}_{+}$-actions on a normal affine surface
$S$ and surjective flat morphisms $q:S\rightarrow C$ with general
fiber isomorphic to $\mathbb{C}$, over nonsingular affine curves
$C$, the latter corresponding to algebraic quotient morphisms associated
with these actions. In this context, Gizatullin \cite{Giz71} and
Bertin \cite{Ber83} (see also \cite{Dub02} for the normal case)
established successively that if a smooth surface $S$ admits an $\mathbb{A}^{1}$-fibration
$q:S\rightarrow C$ as above then this fibration is unique up to isomorphism
of the base if and only if $S$ does not admit a completion $S\hookrightarrow\bar{S}$
by a smooth projective surface $\bar{S}$ for which the boundary divisor
$B=\bar{S}\setminus S$ is zigzag, that is, a chain of nonsingular
rational curves. For instance, the fact that the Danielewski surface
$S_{1}=\left\{ xz-y^{2}+1=0\right\} $ admits two $\mathbb{C}_{+}$-actions
with distinct general orbits can be recovered from this result, as
$S_{1}$ embeds as the complement of a diagonal in $\mathbb{P}^{1}\times\mathbb{P}^{1}$
via the morphism \[
S_{1}\hookrightarrow\mathbb{P}^{1}\times\mathbb{P}^{1},\quad\left(x,y,z\right)\mapsto\left(\left[x:y+1\right],\left[y+1:z\right]\right)=\left(\left[z:y-1\right],\left[x:y-1\right]\right).\]
 Bandman and Makar-Limanov \cite{BML01} (see also \cite{DubG03}
for a more general result) deduced from this criterion that a Danielewski
surface $\rho:S\rightarrow Z_{1,r}$ admits two independent $\mathbb{C}_{+}$-actions
if and only if it is isomorphic to a surface in $\mathbb{C}^{3}$
with equation $xz-P\left(y\right)=0$, where $P$ is a polynomial
with $r$ simple roots. Latter on, Daigle \cite{Dai03} established
that all $\mathbb{C}_{+}$-actions on such a surface $S$ are conjugated
to a one whose general orbits coincides with the ones of the principal
homogeneous $\mathbb{C}_{+}$-bundle structure $\rho:S\rightarrow Z_{1,r}$
factoring the projection $pr_{x}:S\rightarrow\mathbb{C}$. 

\end{enavant}

\begin{enavant} Unfortunately, there is no obvious generalization
of Gizatullin criterion for higher dimensional variety with $\mathbb{C}_{+}$-actions.
However, it turns out that in certain situations such as the one described
in (\ref{thm:Main-Theorem_2}) below, one can establish by direct
computations that the structural morphism $\rho:X\rightarrow Z_{n,r}$
on a Danielewski variety is unique up to automorphisms of the base.
If this holds, then it becomes easier to decide if another Danielewski
variety is isomorphic to $X$ as an abstract variety. Indeed, the
group $\textrm{Aut}\left(Z_{n,r}\right)\times\textrm{Aut}\left(\mathbb{C}_{+}\right)\simeq\textrm{Aut}\left(Z_{n,r}\right)\times\mathbb{C}^{*}$
acts on the set $H^{1}\left(Z_{n,r},\mathcal{O}_{Z_{n,r}}\right)$
by sending a class $\left[g\right]\in H^{1}\left(Z_{n,r},\mathcal{O}_{Z_{n,r}}\right)$
represented by a bundle $\rho:X\rightarrow Z_{n,r}$ with $\mathbb{C}_{+}$-action
$\mu:\mathbb{C}_{+}\times X\rightarrow X$ to the isomorphism class
$\left(\phi,\lambda\right)\cdot\left[g\right]$ of the fiber product
bundle $pr_{2}:\phi^{*}X=X\times_{Z_{n,r}}Z_{n,r}\rightarrow Z_{n,r}$
equipped with the $\mathbb{C}_{+}$-action defined by $\mu_{\lambda}\left(t,\left(x,z\right)\right)\mapsto\left(\mu\left(\lambda^{-1}t,x\right),z\right)$.
Similar arguments as in Theorem 1.1 in \cite{Wil98} imply the following
characterization. 

\end{enavant}

\begin{prop}
\label{pro:Isomorphic-Varieties-Wilkens} Let $\rho_{1}:X_{1}\rightarrow Z_{n,r}$
and $\rho_{2}:X_{2}\rightarrow Z_{n,r}$ be two Danielewski varieties.
If $\rho_{1}$ is a unique $\mathbb{A}^{1}$-bundle structure on $X_{1}$
up to automorphisms of $Z_{n,r}$, then $X_{1}$ and $X_{2}$ are
isomorphic as abstract varieties if their isomorphism classes as principal
$\mathbb{C}_{+}$-bundles belong to the same orbit under the action
of $\textrm{Aut}\left(Z_{n,r}\right)\times\textrm{Aut}\left(\mathbb{C}_{+}\right)$. 
\end{prop}
\begin{enavant} \label{txt:Fibration-on-Dan} Let us again consider
the Danielewski varieties $X_{\left[m\right],\sigma}\subset\mathbb{C}^{n+2}$
with equations $\underline{x}^{\left[m\right]}z-\prod_{i=1}^{r}\left(y-\sigma_{i}\left(\underline{x}\right)\right)=0$,
where $\left[m\right]=\left(m_{1},\ldots,m_{n}\right)\in\mathbb{Z}_{>0}^{n}$
is a multi-index and where $\sigma=\left\{ \sigma_{i}\left(\underline{x}\right)\right\} _{i=1,\ldots,r}$
is collection of polynomials satisfying (1) and (2) in example (\ref{exa:Main-example}).
Again, we denote by $\pi=pr_{\underline{x}}:X_{\left[m\right],\sigma}\rightarrow\mathbb{C}^{n}$,
$\left(\underline{x},y,z\right)\mapsto\underline{x}$ the fibration
which factors through the structural morphism of the principal homogeneous
$\mathbb{C}_{+}$-bundle $\rho:X_{\left[m\right],\sigma}\rightarrow Z_{n,r}$
described in (\ref{exa:Main-example}) above. Suppose that one of
the $m_{i}$'s, say $m_{1}$ is equal to $1$. Then $X_{\left[m\right],\sigma}$
admits a second fibration \[
\pi_{1}:X_{\left[m\right],\sigma}\rightarrow\mathbb{C}^{n},\quad\left(x_{1},\ldots,x_{n},y,z\right)\mapsto\left(x_{2},\ldots,x_{n},z\right)\]
 restricting to the trivial $\mathbb{A}^{1}$-bundle over $\left(\mathbb{C}^{*}\right)^{n}$
and the same argument as in (\ref{exa:Main-example}) above shows
that $\pi_{1}$ factors through the structural morphism of another
principal homogeneous $\mathbb{C}_{+}$-bundle $\rho_{1}:X_{\left[m\right],\sigma}\rightarrow Z_{n,r}$.
On the hand, Makar-Limanov \cite{ML01} established that for every
integer $m\geq2$ the $\mathbb{A}^{1}$-bundle structure $\rho:S\rightarrow Z_{1,r}$
above on a Danielewski surface $S\subset\mathbb{C}^{3}$ with equation
$x^{m}z-P\left(y\right)=0$, where $\deg P\left(y\right)=r\geq2$,
is unique up to isomorphism of the base. More generally, we have the
following result. 

\end{enavant}

\begin{thm}
\label{thm:Main-Theorem_1} Let $\sigma=\left\{ \sigma_{i}\left(\underline{x}\right)\right\} _{i=1,\ldots,r}$
be a collection of $r\geq2$ polynomials satisfying (1) and (2) in
example (\ref{exa:Main-example}). Then for every multi-index $\left[m\right]\in\mathbb{Z}_{>1}^{n}$,
$\rho:X_{\left[m\right],\sigma}\rightarrow Z_{n,r}$ is a unique structure
of principal homogeneous $\mathbb{C}_{+}$-bundle structure on $X_{\left[m\right],\sigma}$
up to action of the group $\textrm{Aut}\left(Z_{n,r}\right)\times\textrm{Aut}\left(\mathbb{C}_{+}\right)$.
\end{thm}
\begin{proof}
This follows from Theorem (\ref{thm:Main-Theorem_2}) below which
guarantees more generally that the algebraic quotient morphism $q:X_{\left[m\right],\sigma}\rightarrow X_{\left[m\right],\sigma}/\!/\mathbb{C}_{+}$
associated with an arbitrary nontrivial $\mathbb{C}_{+}$-action on
$X_{\left[m\right],\sigma}$ coincides with the projection $\pi=pr_{\underline{x}}:X_{\left[m\right],\sigma}\rightarrow\mathbb{C}^{n}$. 
\end{proof}
\noindent It follows from (\ref{txt:Fibration-on-Dan}) that every
Danielewski variety $X_{\left[m\right],\sigma}\subset\mathbb{C}^{n+2}$
defined by a multi-index $\left[m'\right]\in\mathbb{Z}_{\geq1}^{n}\setminus\mathbb{Z}_{>1}^{n}$
admits a second $\mathbb{C}_{+}$-action whose general orbits are
distinct from the general fibers of the $\mathbb{A}^{1}$-bundle $\rho:X_{\left[m\right],\sigma}\rightarrow Z_{n,r}$.
This leads to the following result. 

\begin{cor}
\label{cor:Non-iso_1} For every collection $\sigma=\left\{ \sigma_{i}\left(\underline{x}\right)\right\} _{i=1,\ldots,r}$
of $r\geq2$ polynomials satisfying (1) and (2) in example (\ref{exa:Main-example})
and every pair of multi-index $\left[m\right]\in\mathbb{Z}_{>1}^{n}$
and $\left[m'\right]\in\mathbb{Z}_{\geq1}^{n}\setminus\mathbb{Z}_{>1}^{n}$
the Danielewski varieties $X_{\left[m\right],\sigma}$ and $X_{\left[m'\right],\sigma}$
are not isomorphic.
\end{cor}
\begin{enavant} More generally, let $\left[m\right]=\left(m_{1},\ldots,m_{n}\right)\in\mathbb{Z}_{>1}^{n}$
and $\left[m'\right]=\left(m'_{1},\ldots,m'_{n}\right)\in\mathbb{Z}_{>1}^{n}$
be two multi-index for which the subsets $\left\{ m_{1},\ldots,m_{n}\right\} $
and $\left\{ m'_{1},\ldots,m'_{n}\right\} $ of $\mathbb{Z}$ are
distint. Then for every collection $\sigma=\left\{ \sigma_{i}\left(\underline{x}\right)\right\} _{i=1,\ldots,r}$
of $r\geq2$ polynomials satisfying (1) and (2), the $\check{C}$ech
cocycles \[
g_{ij}=\underline{x}^{-\left[m\right]}\left(\sigma_{j}\left(\underline{x}\right)-\sigma_{i}\left(\underline{x}\right)\right)\quad\textrm{and}\quad g'_{ij}=\underline{x}^{-\left[m'\right]}\left(\sigma_{j}\left(\underline{x}\right)-\sigma_{i}\left(\underline{x}\right)\right)\]
 in $C^{1}\left(\mathcal{U},\mathcal{O}_{Z_{n,r}}\right)\simeq\mathbb{C}\left[\underline{x},\underline{x}^{-1}\right]^{r}$are
not cohomologous and do not belong to the same orbit under the action
of $\textrm{Aut}\left(Z_{n,r}\right)\times\textrm{Aut}\left(\mathbb{C}_{+}\right)$
on $C^{1}\left(\mathcal{U},\mathcal{O}_{Z_{n,r}}\right)$. As a consequence
of Proposition (\ref{pro:Isomorphic-Varieties-Wilkens}) and Theorem
(\ref{thm:Main-Theorem_1}) above, we obtain the following result. 

\end{enavant}

\begin{cor}
\label{cor:Non-Iso_2} Under the hypothesis above, the Danielewski
variety $X_{\left[m\right],\sigma}$ and $X_{\left[m'\right],\sigma}$
are not isomorphic. In particular, there exists an infinite countable
family of pairwise nonisomorphic Danielewski variety $X_{\left[m\right],\sigma}$
with the property that all the varieties $X_{\left[m\right],\sigma}\times\mathbb{C}$
are isomorphic. 
\end{cor}
\begin{rem}
Given a multi-index $\left[m\right]\in\mathbb{Z}_{>1}^{n}$, the problem
of characterizing explicitly the collections $\sigma=\left\{ \sigma_{i}\left(\underline{x}\right)\right\} _{i=1,\ldots,r}$
which lead to isomorphic Danielewski varieties $X_{\left[m\right],\sigma}$
is more subtle in general. By virtue of proposition (\ref{pro:Isomorphic-Varieties-Wilkens}),
it is equivalent to describe the orbits of the associated cocycles
$g_{ij}=\underline{x}^{-\left[m\right]}\left(\sigma_{j}\left(\underline{x}\right)-\sigma_{i}\left(\underline{x}\right)\right)$
under the action of $\textrm{Aut}\left(Z_{n,r}\right)\times\textrm{Aut}\left(\mathbb{C}_{+}\right)$.
In the case of surfaces, the question becomes simpler as $\textrm{Aut}\left(Z_{1,r}\right)\simeq\mathbb{C}^{*}\times\mathbb{Z}/r\mathbb{Z}$.
For instance, Makar-Limanov \cite{ML01} obtained a complete classification
of the Danielewski surfaces $S\subset\mathbb{C}^{3}$ with equation
$x^{n}z-P\left(y\right)=0$, where $n\geq2$. More generally, we refer
the interested reader to the forthcoming paper \cite{DubP04}, in
which we study Danielewski surfaces with equations $x^{n}z-Q\left(x,y\right)=0$.
\end{rem}

\section{Additive group actions on Danielewski varieties}

Makar-Limanov \cite{ML96} observed that it is sometimes possible
to obtain information on algebraic $\mathbb{C}_{+}$-actions on an
affine variety $X$ by considering homogeneous $\mathbb{C}_{+}$-actions
on certain affine cones $\hat{X}$ associated with $X$. We recall
that \emph{the Makar-Limanov invariant} of an affine variety $X=\textrm{Spec}\left(B\right)$
is the subring $\textrm{ML}\left(X\right)$ of $B$ consisting of
regular functions on $B$ which are invariant under all $\mathbb{C}_{+}$-actions
on $X$. Using associated homogeneous objects, he established in \emph{}\cite{ML96}
\emph{}that the Makar-Limanov invariant of the Russell cubic threefold,
\emph{i.e.} the hypersurface $X\subset\mathbb{C}^{4}$ with equation
$x+x^{2}y+z^{2}+t^{3}=0$, is not trivial. He also computed in \cite{ML01}
the Makar-Limanov invariant of affine surfaces $S=\left\{ x^{n}z-P\left(y\right)=0\right\} $,
where $\deg\left(P\right)>1$and $n>1$. Here we use a similar method,
based on real-valued weight degree functions, to compute the Makar-Limanov
invariant of the Danielewski varieties $X_{\left[m\right],\sigma}$,
where $\left[m\right]\in\mathbb{Z}_{>1}^{n}$.

\subsection{Basic facts on locally nilpotent derivations}

\indent\newline\indent He we recall results on locally nilpotent
derivations that will be used in the following subsections. We refer
the reader to \cite{FLN92} and \cite{KML97} for more complete discussions. 

\begin{enavant} Algebraic $\mathbb{C}_{+}$-actions on a complex
affine variety $X=\textrm{Spec}\left(B\right)$ are in one-to-one
correspondence with locally nilpotent $\mathbb{C}$-derivations of
$B$, that is, derivations $\partial:B\rightarrow B$ such that every
$b$ belongs to the kernel of $\partial^{m}$ for a suitable $m=m\left(b\right)$.
Indeed, for every algebraic $\mathbb{C}_{+}$-action on $S$ with
comorphism $\mu^{*}:B\rightarrow B\otimes_{\mathbb{C}}\mathbb{C}\left[t\right]$,
$\partial_{\mu}={\displaystyle \frac{d}{dt}\mid_{t=0}\circ\mu^{*}:B\rightarrow B}$
is a locally nilpotent derivation. Conversely, for every such derivation
$\partial:B\rightarrow B$ the exponential map \[
\textrm{exp}\left(t\partial\right):B\rightarrow B\left[t\right],\quad b\mapsto\sum_{n\geq0}\frac{\partial^{n}b}{n!}t^{n}\]
 coincides with the comorphism of an algebraic $\mathbb{C}_{+}$-action
on $X$. To every locally nilpotent derivation $\partial$ of $B$,
we associate a function \[
\deg_{\partial}:B\rightarrow\mathbb{N}\cup\left\{ -\infty\right\} ,\textrm{ defined by }\deg_{\partial}\left(b\right)=\begin{cases}
-\infty & \textrm{if }b=0\\
\max\left\{ m,\,\partial^{m}b\neq0\right\}  & \textrm{otherwise},\end{cases}\]
 which we call the \emph{degree function generated by} $\partial$.
We recall the following facts. 

\end{enavant} 

\begin{prop}
\label{pro:LND-facts} Let $\partial$ be a nontrivial locally nilpotent
derivation of $B$. Then the following hold.

(1) $B$ has transcendence degree one over $\textrm{Ker}\left(\partial\right)$.
The field of fraction $Frac\left(B\right)$ of $B$ is a purely transcendental
extension of $Frac\left(\textrm{Ker}\left(\partial\right)\right)$,
and $\textrm{Ker}\left(\partial\right)$ is algebraically closed in
$B$.

(2) For every $f\in\textrm{Ker}\left(\partial^{2}\right)\setminus\textrm{Ker}\left(\partial\right)$,
the localization $B_{f}$ of $B$ at $f$ is isomorphic to the polynomial
ring in one variable $\textrm{Ker}\left(\partial\right)_{\partial\left(f\right)}\left[f\right]$
over the localization $\textrm{Ker}\left(\partial\right)_{\partial\left(f\right)}$
of $\textrm{Ker}\left(\partial\right)$ at $\partial\left(f\right)$.
In particular, for every $b\in\textrm{Ker}\left(\partial^{m+1}\right)\setminus\textrm{Ker}\left(\partial^{m}\right)$,
there exists $a',a_{0},\ldots,a_{m}\in\textrm{Ker}\left(\partial\right)$,
where $a',a_{m}\neq0$, such that $a'b=\sum_{j=0}^{m}a_{j}f^{j}$. 

(3) $\deg_{\partial}:B\rightarrow\mathbb{N}\cup\left\{ -\infty\right\} $
is a degree function, i.e. $\deg_{\partial}\left(b+b'\right)\leq\max\left(\deg_{\partial}\left(b\right),\deg_{\partial}\left(b'\right)\right)$
and $\deg_{\partial}\left(bb'\right)=\deg_{\partial}\left(b\right)+\deg_{\partial}\left(b'\right)$. 

(4) If $b,b'\in B\setminus\left\{ 0\right\} $ and $bb'\in\textrm{Ker}\left(\partial\right)$,
then $b,b'\in\textrm{Ker}\left(\partial\right)$.
\end{prop}

\subsection{Equivariant deformations to the cone following Kaliman and Makar-Limanov}

\indent\newline\indent Here we review a procedure due to Makar-Limanov
\cite{ML96} which associates to filtered algebra $\left(B,\mathcal{F}\right)$
equipped with a locally nilpotent derivation $\partial$ a graded
algebra equipped with an homogeneous locally nilpotent derivation
induced by $\partial$.

\begin{enavant} We let $B$ be a finitely generated algebra, equipped
with an exhaustive, separated, ascending filtration $\mathcal{F}=\left\{ F^{t}B\right\} _{t\in\mathbb{R}}$
by $\mathbb{C}$-linear subspaces $F^{t}B$ of $B$. For every $t\in\mathbb{R}$,
we let $F_{0}^{t}B=\bigcup_{s<t}F^{s}B$. We denote by \begin{eqnarray*}
gr_{\mathcal{F}}B & = & \bigoplus_{t\in\mathbb{R}}\left(gr_{\mathcal{F}}B\right)^{t},\quad\textrm{where }\left(gr_{\mathcal{F}}B\right)^{t}=F^{t}B/F_{0}^{t}B\end{eqnarray*}
 the $\mathbb{R}$-graded algebra associated to the filtered algebra
$\left(B,\mathcal{F}\right)$, and we let $gr:B\rightarrow gr_{\mathcal{F}}B$
the natural map which sends an element $b\in F^{t}B\subset B$ to
its image $gr\left(b\right)$ under the canonical map $F^{t}B\rightarrow F^{t}B/F_{0}^{t}B\subset gr_{\mathcal{F}}B$.
Suppose further that $1\in F^{0}B\setminus F_{0}^{0}B$ and that \[
\left(F^{t_{1}}B\setminus F_{0}^{t_{1}}B\right)\left(F^{t_{2}}B\setminus F_{0}^{t_{2}}B\right)\subset\left(F^{t_{1}+t_{2}}B\setminus F_{0}^{t_{1}+t_{2}}B\right)\qquad\textrm{for every }t_{1},t_{2}\in\mathbb{R}.\]
 Then the filtration $\mathcal{F}$ is induced by a degree function
$d_{\mathcal{F}}:B\rightarrow\mathbb{R}\cup\left\{ -\infty\right\} $
on $B$. Indeed, the formulas $d_{\mathcal{F}}\left(0\right)=-\infty$
and $d_{\mathcal{F}}\left(b\right)=t$ if $b\in F^{t}B\setminus F_{0}^{t}B\subset B$
define a degree function on $B$ such that $F^{t}B=\left\{ b\in B,\, d\left(b\right)\leq t\right\} $
for every $t\in\mathbb{R}$. In what follows, we only consider filtrations
induced by degree functions. 

\end{enavant}

\begin{enavant} Given a nontrivial locally nilpotent derivation $\partial$
of $B$ and a nonzero $b\in B$, we let $t\left(b\right)=d_{\mathcal{F}}\left(\partial b\right)-d_{\mathcal{F}}\left(b\right)\in\mathbb{R}$.
By definition, if $b\in F^{t}B\setminus\left(\textrm{Ker}\partial\cap F_{0}^{t}B\right)$
then $\partial b\in F^{t+t\left(b\right)}B\setminus F_{0}^{t+t\left(b\right)}B$.
Since $B$ is finitely generated, it follows that there exists a smallest
$t_{0}\in\mathbb{R}$ such that $\partial F^{t}B\subset F^{t+t_{0}}B$.
So $\partial$ induces a locally nilpotent derivation $gr\partial$
of the associated graded algebra $gr_{\mathcal{F}}B$ of $\left(B,\mathcal{F}\right)$,
defined by \begin{eqnarray*}
gr\partial\left(gr\left(b\right)\right) & = & \begin{cases}
gr\left(\partial b\right) & \textrm{if }d_{\mathcal{F}}\left(\partial\left(b\right)\right)-d_{\mathcal{F}}\left(b\right)=t_{0}\\
0 & \textrm{otherwise}.\end{cases}\end{eqnarray*}
 By construction, $gr\partial$ sends an homogeneous component $F^{t}B/F_{0}^{t}B$
of $gr_{\mathcal{F}}B$ into the homogeneous component $F^{t+t_{0}}B/F_{0}^{t+t_{0}}B$.
We say that $gr\partial$ is the \emph{homogeneous locally nilpotent
derivation of} $gr_{\mathcal{F}}B$ \emph{associated with} $\partial$.
By construction, if $gr_{\mathcal{F}}B$ is a domain, then \begin{eqnarray}
\deg_{\partial}\left(b\right) & \geq & \deg_{gr\partial}\left(gr\left(b\right)\right)\label{eq:degree-comp}\end{eqnarray}
 for every $b\in B$. We will see below that this inequality plays
a crucial role in the computation of the Makar-Limanov invariant of
certain Danielewski varieties. 

\end{enavant}

\begin{rem}
For integral-valued degree functions $d:B\rightarrow\mathbb{Z}\cup\left\{ -\infty\right\} $,
the above construction admits a simple geometric interpretation. Indeed,
letting $\mathcal{F}=\left\{ F^{n}B\right\} _{n\in\mathbb{Z}}$ be
the filtration generated by $d$, we consider the Rees algebra \[
\mathcal{R}\left(B,\mathcal{F}\right)=\bigoplus_{n\in\mathbb{Z}}\mathcal{F}^{n}s^{-n}\subset B\left[s,s^{-1}\right].\]
 Every locally nilpotent derivation $\partial$ of $B$ canonically
extends to a locally nilpotent derivation $\tilde{\partial}$ of $\mathcal{R}\left(B,\mathcal{F}\right)$
with the property that $\tilde{\partial}\left(s\right)=0$. By construction,
the inclusion $\mathbb{C}\left[s\right]\hookrightarrow\mathcal{R}\left(B,\mathcal{F}\right)$
gives rise to a flat family $\rho:\mathcal{X}=\textrm{Spec}\left(\mathcal{R}\left(B,\mathcal{F}\right)\right)\rightarrow\mathbb{C}$
of affine varieties with $\mathbb{C}_{+}$-actions, such that for
every $s\in\mathbb{C}^{*}$, the fiber $\mathcal{X}_{s}$ is isomorphic
to $X$ equipped with the $\mathbb{C}_{+}$-action defined by $\partial$,
whereas the fiber $\mathcal{X}_{0}\simeq\textrm{Spec}\left(\mathcal{R}\left(B,\mathcal{F}\right)/s\mathcal{R}\left(B,\mathcal{F}\right)\right)$
is isomorphic to the spectrum of the graded algebra $gr_{\mathcal{F}}B$,
equipped with $\mathbb{C}_{+}$-action corresponding to the homogeneous
locally nilpotent derivation $gr\partial$ of $gr_{\mathcal{F}}B$
defined above. 
\end{rem}

\subsection{On the Makar-Limanov invariants of Danielewski varieties $X_{\left[m\right],\sigma}$ }

\indent\newline\indent Here we consider a class of affine varieties
with $\mathbb{C}_{+}$-actions which contains the Danielewski varieties
$X_{\left[m\right],\sigma}$ of example (\ref{exa:Main-example}).
We construct certain filtrations $\mathcal{F}_{d}$ of their coordinate
rings induced by weight degree functions, and we determine the structure
of the associated homogeneous objects. Finally we compute their Makar-Limanov
invariants. 

\begin{defn}
\label{def:Var-def} Given a monic polynomial $Q\left(\underline{x},y\right)=y^{r}+\sum_{i=0}^{r-1}a_{i}\left(\underline{x}\right)y^{i}\in\mathbb{C}\left[\underline{x}\right]\left[y\right]$
of degree $r\geq2$ and a multi-index $\left[m\right]=\left(m_{1},\ldots,m_{n}\right)\in\mathbb{Z}_{\geq1}^{n}$,
we denote by $X_{\left[m\right],Q}\subset\mathbb{C}^{n+2}$ the affine
variety with equation $\underline{x}^{\left[m\right]}z-Q\left(\underline{x},y\right)=0$. 
\end{defn}
\begin{enavant} \label{txt:Fibration-on-Vars} Clearly, the above
class of affine varieties contains the Danielewski varieties $X_{\left[m\right],\sigma}\subset\mathbb{C}^{n+2}$
with equations $\underline{x}^{\left[m\right]}z-\prod_{i=1}^{r}\left(y-\sigma_{i}\left(\underline{x}\right)\right)=0$.
Again, the projection \[
\pi=pr_{\underline{x}}:X_{\left[m\right],Q}\rightarrow\mathbb{C}^{n},\quad\left(\underline{x},y,z\right)\mapsto\underline{x}\]
 restricts to a trivial $\mathbb{A}^{1}$-bundle $\left(\mathbb{C}^{*}\right)^{n}\times\mathbb{C}=\textrm{Spec}\left(\mathbb{C}\left[\underline{x},\underline{x}^{-1}\right]\left[y\right]\right)$
over $\left(\mathbb{C}^{*}\right)^{n}\subset\mathbb{C}^{n}$. On the
other hand $\pi^{-1}\left(H_{\underline{x}}\right)_{red}$ is the
disjoint union of $\tilde{r}$ copies of $H_{\underline{x}}\times\mathbb{C}$
with equations $\left\{ x_{1}\cdots x_{n}=0,\, y=\lambda_{i}\right\} $,
where $\lambda_{1},\ldots,\lambda_{\tilde{r}}$ denote the distinct
roots of the polynomial $P\left(y\right)=Q\left(0,y\right)$. The
locally nilpotent derivation $\partial$ of $\mathbb{C}\left[\underline{x},y,z\right]$
defined by \[
\partial\left(x_{i}\right)=0,\, i=1,\ldots,r,\quad\partial\left(y\right)=\underline{x}^{\left[m\right]}\quad\partial\left(z\right)=\frac{\partial Q\left(x,y\right)}{\partial y}\]
 annihilates the definning ideal $I=\left(\underline{x}^{\left[m\right]}z-Q\left(x,y\right)\right)$
of $X_{\left[m\right],Q}$, whence induces a nontrivial locally nilpotent
derivation of the coordinate ring $B$ of $X_{\left[m\right],Q}$.
The general orbits of the corresponding $\mathbb{C}_{+}$-action coincide
with the general fibers of $\pi$. Hence $\pi$ coincides with the
algebraic quotient morphism $q:X_{\left[m\right],Q}\rightarrow X_{\left[m\right],Q}/\!/\mathbb{C}_{+}=\textrm{Spec}\left(B^{\mathbb{C}_{+}}\right)$.
This shows that $\textrm{ML}\left(X_{\left[m\right],Q}\right)\subset\mathbb{C}\left[\underline{x}\right]$.
Actually, a similar argument as in (\ref{txt:Fibration-on-Dan}) above
shows that $\textrm{ML}\left(X_{\left[m\right],Q}\right)$ is a subring
of $\mathbb{C}\left[x_{i_{1}},\ldots,x_{i_{s}}\right]$, where $i_{1},\ldots,i_{s}$
denote the indices for which $m_{i_{k}}=1$. In particular, if $\left[m\right]=\left(1,\ldots,1\right)$,
then $\textrm{ML}\left(X_{\left[m\right],Q}\right)=\mathbb{C}$. In
contrast, we have the following result. 

\end{enavant}

\begin{thm}
\label{thm:Main-Theorem_2} If $\left[m\right]\in\mathbb{Z}_{>1}^{n}$
then the Makar-Limanov invariant of a variety $X_{\left[m\right],Q}$
is isomorphic to $\mathbb{C}\left[\underline{x}\right]$. 
\end{thm}
\begin{enavant} It suffices to shows $\textrm{Ker}\left(\partial^{2}\right)\subset\mathbb{C}\left[\underline{x},y\right]\subset B$
for every nontrivial locally nilpotent derivation $\partial$ on the
coordinate ring $B$ of $X_{\left[m\right],Q}$. Indeed, if $\partial$
is nontrivial, then it follows from (2) in Proposition (\ref{pro:LND-facts})
that there exists $f\in\textrm{Ker}\left(\partial^{2}\right)\setminus\textrm{Ker}\left(\partial\right)$
such that $z=\underline{x}^{-\left[m\right]}\left(y^{2}-1\right)\in B\subset\mathbb{C}\left[\underline{x},\underline{x}^{-1},y\right]$
satisfies a relation of the form $a'z=\sum_{j=1}^{m}a_{j}f^{j}$ for
suitable elements $a',a_{0},\ldots,a_{m}\in\textrm{Ker}\left(\partial\right)$,
where $a',a_{m}\neq0$. Therefore, if $\textrm{Ker}\left(\partial^{2}\right)\subset\mathbb{C}\left[\underline{x},y\right]$
then $z=r\left(\underline{x},y\right)/q\left(\underline{x},y\right)$
for a certain polynomial $q\left(\underline{x},y\right)\in\textrm{Ker}\left(\partial\right)$.
This implies that $\underline{x}^{\left[m\right]}$ divides $q\left(\underline{x},y\right)$
and so, by virtue of (3) in (\ref{pro:LND-facts}), $\mathbb{C}\left[\underline{x}\right]\subset\textrm{Ker}\left(\partial\right)$
as $m_{i}\geq1$ for every $i=1,\ldots,n$. To show that the inclusion
$\textrm{Ker}\left(\partial^{2}\right)\subset\mathbb{C}\left[\underline{x},y\right]$
holds for every nontrivial locally nilpotent derivation on $B$, we
study in (\ref{def:weight-degree})-(\ref{cor:Var-kernels}) below
the homogeneous objects associated with certain filtrations on $B$
induced by weight degree functions.

\end{enavant}

\begin{defn}
\label{def:weight-degree} A \emph{weight degree function} on a polynomial
ring $\mathbb{C}\left[\underline{x}\right]$ is a degree function
$d:\mathbb{C}\left[\underline{x}\right]\rightarrow\mathbb{R}$ defined
by real weights $d_{i}=d\left(x_{i}\right)$, $i=1,\ldots,n$. The
$d$-degree of monomial $m=x^{\left[\alpha\right]}$ is $\alpha_{1}d_{1}+\ldots+\alpha d_{n}$,
and the $d$-degree $d\left(p\right)$ of a polynomial $p\in\mathbb{C}\left[\underline{x}\right]$
is defined as the suppremum of the degrees $d\left(m\right)$, where
$m$ runs through the monomials of $p$. A weight degree function
$d$ defines a grading $\mathbb{C}\left[\underline{x}\right]=\bigoplus_{t\in\mathbb{R}}\mathbb{C}\left[\underline{x}\right]_{t}$,
where $\mathbb{C}\left[\underline{x}\right]_{t}\setminus\left\{ 0\right\} $
consists of all the $d$-homogeneous polynomials of $d$-degree $t$.
In what follows, we denote by $\bar{p}$ the \emph{principal} $d$-\emph{homogeneous
component of} $p$, that is, the homogeneous component of $p$ of
degree $d\left(p\right)$. A degree function $d$ on $\mathbb{C}\left[\underline{x}\right]$
naturally extends to a degree function on the algebra $\mathbb{C}\left[\underline{x},\underline{x}^{-1}\right]$
of Laurent polynomials. 
\end{defn}
\begin{enavant} Given a multi-index $\left[m\right]\in\mathbb{Z}_{>1}^{n}$
and a monic polynomial $Q\left(\underline{x},y\right)\in\mathbb{C}\left[\underline{x}\right]\left[y\right]$
as in definition (\ref{def:Var-def}), we denote by $B=\mathbb{C}\left[\underline{x},y,z\right]/I$,
where $I=\left(\underline{x}^{\left[r\right]}z-Q\left(\underline{x},y\right)\right)$,
the coordinate ring of the corresponding variety $X_{\left[m\right],Q}$,
and we denote by $\sigma:\mathbb{C}\left[\underline{x},y,z\right]\rightarrow B$
the natural morphism. The polynomial ring $\mathbb{C}\left[\underline{x},y\right]$
is naturally a subring of $B$. Moreover, by means of the localization
homomorphism $B\hookrightarrow B_{\underline{x}}=B\otimes_{\mathbb{C}\left[\underline{x}\right]}\mathbb{C}\left[\underline{x},\underline{x}^{-1}\right]\simeq\mathbb{C}\left[\underline{x},\underline{x}^{-1},y\right]$,
$B$ is itself identified to the subalgebra $\mathbb{C}\left[\underline{x},y,\underline{x}^{-\left[m\right]}Q\left(\underline{x},y\right)\right]$
of $\mathbb{C}\left[\underline{x},\underline{x}^{-1},y\right]$. Hence
every weight degree function $d$ on $\mathbb{C}\left[\underline{x},\underline{x}^{-1},y\right]$
induces an exhaustive separated ascending filtration $\mathcal{F}_{d}=\left\{ F^{t}B\right\} _{t\in\mathbb{R}}$
of $B\subset\mathbb{C}\left[\underline{x},\underline{x}^{-1},y\right]$
by means of the subsets $F^{t}B=\left\{ p\in B,\: d\left(p\right)\leq t\right\} $,
$t\in\mathbb{R}$. 

\end{enavant}

\begin{enavant} \label{txt:Graded-Iso} Since $Q\left(\underline{x},y\right)=y^{r}+\sum_{i=0}^{r-1}a_{i}\left(\underline{x}\right)y^{i}$
is monic, it follows that if the weight $d_{y}$ of $y$ is positive
and sufficiently bigger that the weights $d_{i}$ of the $x_{i}$'s,
then the principal $d$-homogeneous component of $Q\left(\underline{x},y\right)$
is simply $\bar{Q}\left(\underline{x},y\right)=y^{r}$. If this holds,
then $gr_{\mathcal{F}}B$ is generated by $gr\left(x\right)=x$, $gr\left(y\right)=y$
and $gr\left(z\right)=x^{-\left[m\right]}y^{r}$, with the unique
relation $x^{\left[m\right]}gr\left(z\right)=y^{r}$. Hence, letting
$\tilde{d}:\mathbb{C}\left[\underline{x},y,z\right]\rightarrow\mathbb{R}$
be the unique weight degree function restricting to $d$ on $\mathbb{C}\left[\underline{x},y\right]\subset\mathbb{C}\left[\underline{x},y,z\right]$
and such that $\tilde{d}\left(z\right)=rd_{y}-\left(m_{1}d_{1}+\dots+m_{n}d_{n}\right)\in\mathbb{R}$,
we obtain an isomorphism of graded algebras \[
\phi:\hat{B}=\mathbb{C}\left[\underline{x},y,z\right]/\hat{I}=\bigoplus_{t\in\mathbb{R}}\hat{B}^{t}\stackrel{\sim}{\longrightarrow}gr_{\mathcal{F}_{d}}B=\bigoplus_{t\in\mathbb{R}}F^{t}B/F_{0}^{t}B,\]
 where $\hat{I}=\left(\underline{x}^{\left[m\right]}z-y^{r}\right)\subset\mathbb{C}\left[\underline{x},y,z\right]$
denotes the $\tilde{d}$-homogeneous ideal generated by the principal
components of the polynomials in $I=\left(\underline{x}^{\left[r\right]}z-Q\left(\underline{x},y\right)\right)$,
and where $\hat{B}^{t}=\hat{B}^{t}=\mathbb{C}\left[\underline{x},y,z\right]_{t}/\hat{I}\cap\mathbb{C}\left[\underline{x},y,z\right]_{t}$
for every $t\in\mathbb{R}$. 

\end{enavant}

\begin{enavant} \label{txt:Special-weights} It follows from (1)
and (4) in Proposition (\ref{pro:LND-facts}), that the kernel of
an associated homogeneous locally nilpotent derivations $gr\partial$
of $gr_{\mathcal{F}_{d}}B$ contains $n$ algebraically independent
irreducible homogeneous elements. To make the study of these derivations
easier, we need to make the set of these irreducible homogeneous elements
as small as possible. For this purpose, we consider weight functions
$d:\mathbb{C}\left[\underline{x},y\right]\rightarrow\mathbb{R}$ satisfying
the following properties :\\

(1) The weight $d_{y}$ of $y$ is positive, and $\bar{Q}\left(\underline{x},y\right)=y^{r}$. 

(2) The real weights $d_{i}=d\left(x_{i}\right)$ and $d_{y}$ are
linearly independent over $\mathbb{Z}$. \\

\noindent According to (\ref{txt:Graded-Iso}) above, the first condition
guarantees that the graded algebra $gr_{\mathcal{F}}B$ of the filtered
algebra $\left(B,\mathcal{F}_{d}\right)$ is isomorphic to the quotient
$\hat{B}$ of $\mathbb{C}\left[\underline{x},y,z\right]$ by the $\tilde{d}$-homogeneous
ideal $\hat{I}=\left(\underline{x}^{\left[m\right]}z-y^{r}\right)$.
The second one is motivated by the following result. 

\end{enavant} 

\begin{lem}
\label{lem:Variables} Under the hypothesis above, every homogeneous
element of $\hat{B}$ is the image by the natural morphism $\hat{\sigma}:\mathbb{C}\left[\underline{x},y,z\right]\rightarrow\hat{B}$
of a unique monomial of $\mathbb{C}\left[\underline{x},y,z\right]$
not divisible by $\underline{x}^{\left[m\right]}z$. In particular,
every irreducible homogeneous element of $\hat{B}$ is the image of
a variable of $\mathbb{C}\left[\underline{x},y,z\right]$. 
\end{lem}
\begin{proof}
Since $\hat{I}=\left(\underline{x}^{\left[m\right]}z-y^{r}\right)$,
every nonzero homogeneous element of $\hat{B}$ is the image by $\hat{\sigma}$
of a unique homogeneous polynomial $\bar{p}\in\mathbb{C}\left[\underline{x},y,z\right]$
whose monomials are not divisible by $\underline{x}^{\left[m\right]}z$.
On the other hand, the hypothesis on $d,$ together with the fact
that $\tilde{d}\left(z\right)=2d_{y}-\left(m_{1}d_{1}+\ldots+m_{n}d_{n}\right)$
implies that if $\bar{p}$ contains a pair of monomials $\mu_{1}\neq\mu_{2}$,
then there exists $\lambda\in\mathbb{C}$ and $k\in\mathbb{Z}$ such
that $\mu_{1}\mu_{2}^{-1}=\lambda\left(\underline{x}^{\left[m\right]}zy^{-r}\right)^{k}$.
If $k\neq0$, then $\underline{x}^{\left[m\right]}z$ divides one
of the $\mu_{i}$, which is impossible. Thus $\bar{p}$ is a monomial. 
\end{proof}
\begin{prop}
\label{pro:Homogeneous-kernel} If $\left[m\right]\in\mathbb{Z}_{>1}^{n}$
then $\textrm{Ker}\left(\hat{\partial}\right)=\mathbb{C}\left[\underline{x}\right]$
for every associated homogeneous locally nilpotent derivation $\hat{\partial}$
on $\hat{B}$. Furthermore $\deg_{\hat{\partial}}\left(\hat{\sigma}\left(z\right)\right)\geq2$. 
\end{prop}
\begin{proof}
By virtue of (1) and (4) in (\ref{pro:LND-facts}), the kernel of
$\hat{\partial}$ contains $n$ algebraically independent irreducible
homogeneous elements $\xi_{1},\ldots,\xi_{n}$. So it follows from
lemma (\ref{lem:Variables}) above that the $\xi_{i}$'s are the images
by $\hat{\sigma}$ of $n$ distinct variables of $\mathbb{C}\left[\underline{x},y,z\right]$.
These functions $\xi_{i}$, $i=1,\ldots,n$, define a morphism $q:\hat{X}=\textrm{Spec}\left(\hat{B}\right)\rightarrow\mathbb{C}^{n}$
which invariant for the $\mathbb{C}_{+}$-action defined by $\hat{\partial}$.
In particular, for a general point $\lambda=\left(\lambda_{1},\ldots,\lambda_{n}\right)\in\mathbb{C}^{n}$,
the $\mathbb{C}_{+}$-action on $\hat{X}$ specializes to a nontrivial
$\mathbb{C}_{+}$-action on the fiber $q^{-1}\left(\lambda\right)$.
Suppose that one of the $\xi_{i}$'s, say $\xi_{1}$, is the image
of $y$. Then, depending on the other variables inducing the $\xi_{i}$'s,
$i=2,\ldots,n$, we would obtain, for a general $\mu\in\mathbb{C}$,
a nontrivial $\mathbb{C}_{+}$-action on one of the curves $C\subset\mathbb{C}^{2}$
with equations $x_{i_{1}}^{m_{i_{1}}}x_{i_{2}}^{m_{i_{2}}}-\mu=0$
or $x_{i_{1}}^{m_{i_{1}}}z-\mu=0$, which is absurd. Similarly, is
$\xi_{1}$ is the image of $z$ then, for a general $\mu\in\mathbb{C}$,
the $\mathbb{C}_{+}$-action on $\hat{X}$ would specialize to a nontrivial
action on the curve with equation $\mu x_{i}^{m_{i}}-y^{r}=0$ for
certain $i=1,\ldots,n$. This impossible as $r>1$ and $m_{i}>1$
for every $i=1,\ldots,n$ by hypothesis. This proves that $\textrm{Ker}\left(\hat{\partial}\right)$
contains $\mathbb{C}\left[\underline{x}\right]$. Thus $\hat{\partial}$
naturally extends to a locally nilpotent derivation of $\hat{B}_{\underline{x}}\simeq\mathbb{C}\left[\underline{x},\underline{x}^{-1},y\right]$.
In turn, this implies that $\deg_{\hat{\partial}}\left(y\right)=1$
and $\deg_{\hat{\partial}}\left(\hat{\sigma}\left(z\right)\right)\geq2$
as $\hat{\sigma}\left(z\right)\in\hat{B}$ coincides with $x^{-\left[m\right]}y^{r}\in\hat{B}_{\underline{x}}$
via the canonical injection $\hat{B}\hookrightarrow\hat{B}_{\underline{x}}$.
Therefore, the projection $pr_{\underline{x}}:\hat{X}\rightarrow\mathbb{C}^{n}$
coincides with the algebraic quotient morphism of the associated $\mathbb{C}_{+}$-action.
This proves that $\textrm{Ker}\left(\hat{\partial}\right)=\mathbb{C}\left[\underline{x}\right]$. 
\end{proof}
\noindent The following result completes the proof of Theorem (\ref{thm:Main-Theorem_2}).

\begin{cor}
\label{cor:Var-kernels} For every nontrivial locally nilpotent $\partial$
of $B$, $\textrm{Ker}\left(\partial^{2}\right)$ is contained in
$\mathbb{C}\left[\underline{x},y\right]$. 
\end{cor}
\begin{proof}
Recall that $b\in\textrm{Ker}\left(\partial^{2}\right)$ if and only
if $\deg_{\partial}\left(b\right)\leq1$. Since $I$ is generated
by the polynomial $\underline{x}^{\left[m\right]}z-Q\left(\underline{x},y\right)$,
every $b\in\textrm{Ker}\left(\partial^{2}\right)$ is the restriction
to $X_{\left[m\right],Q}$ of a unique polynomial $p\in\mathbb{C}\left[\underline{x},y,z\right]$
whose monomials are not divisible by $\underline{x}^{\left[m\right]}z$.
Suppose that $p\not\in\mathbb{C}\left[\underline{x},y\right]$. Then
there exists a weight degree function $d$ on $\mathbb{C}\left[\underline{x},y,z\right]$
as in (\ref{txt:Special-weights}) for which the principal $d$-homogeneous
component $\bar{p}$ belongs to $\mathbb{C}\left[\underline{x},y,z\right]\setminus\mathbb{C}\left[\underline{x},y\right]$.
We deduce from lemma (\ref{lem:Variables}) above that $\bar{p}=\underline{x}^{\left[\alpha\right]}y^{\beta}z^{\gamma}$,
where $\gamma\geq1$ and $\underline{x}^{\left[m\right]}z$ does not
divide $\underline{x}^{\left[\alpha\right]}z^{\gamma}$. Letting $\hat{\partial}=gr\partial$
be the homogeneous locally nilpotent derivation of $\hat{B}=gr_{\mathcal{F}}B$
associated with $\partial$, we have $\deg_{\hat{\partial}}\left(\hat{\sigma}\left(\bar{p}\right)\right)\geq\deg_{\hat{\partial}}\left(\hat{\sigma}\left(z\right)\right)$
and so (see (\ref{eq:degree-comp})), $\deg_{\partial}\left(b\right)\geq\deg_{\hat{\partial}}\left(\hat{\sigma}\left(z\right)\right)$
as $\hat{\sigma}\left(\bar{p}\right)$ coincides via the isomorphism
$\phi$ of (\ref{txt:Graded-Iso}) with the image $gr\left(b\right)\in gr_{\mathcal{F}}B$
of $b$. This is absurd as $\deg_{\hat{\partial}}\left(\hat{\sigma}\left(z\right)\right)\geq2$
by virtue of lemma (\ref{pro:Homogeneous-kernel}). 
\end{proof}
\bibliographystyle{amsplain}

\begin{thebibliography}{10}

\bibitem{BML01}
T.~Bandman and L.~Makar-Limanov, \emph{Affine surfaces with
  ${AK}\left({S}\right)=\mathbb{C}$}, Michigan J. Math. \textbf{49} (2001),
  567--582.

\bibitem{Ber83}
J.~Bertin, \emph{Pinceaux de droites et automorphismes des surfaces affines},
  J. reine angew. Math. \textbf{341} (1983), 32--53.

\bibitem{Dai03}
D.~Daigle, \emph{On locally nilpotent derivations of $k\left[x,y,z\right] /
  \left(xy-p\left(z\right)\right)$}, JPPA \textbf{181} (2003), 181--208.

\bibitem{Dan89}
W.~Danielewski, \emph{On a cancellation problem and automorphism groups of
  affine algebraic varieties}, Preprint {W}arsaw, 1989.

\bibitem{Dub02}
A.~Dubouloz, \emph{Completions of normal affine surfaces with a trivial
  {M}akar-{L}imanov invariant}, Michigan J. Math. \textbf{52} (2004), no.~2,
  289--308.

\bibitem{DubG03}
A.~Dubouloz, \emph{{D}anielewski-{F}ieseler {S}urfaces}, {T}ransformation {G}roups
  \textbf{10} (2005), no.~2, 139--162.

\bibitem{DubP04}
A.~Dubouloz and P.M. Poloni, \emph{Automorphisms of {D}anielewski
  hypersurfaces}, In preparation.

\bibitem{EaH73}
P.~Eaking and W.~Heinzer, \emph{A cancellation theorem for rings}, Conference
  on {C}ommutative {A}lgebra. Lecture {N}otes in {M}athematics vol. 311 (J.W.
  Brewer and E.A. Rutter, eds.), {S}pringer {V}erlag,
  {B}erlin-{H}eidelberg-{N}ew {Y}ork, 1973, pp.~61--77.

\bibitem{FLN92}
Y.~Ferrero, M.~Lequain and A.~Nowicki, \emph{A note on locally nilpotent
  derivations}, J. of {P}ure and {A}ppl {A}lgebra \textbf{79} (1992), 45--50.

\bibitem{Fie94}
K.H. Fieseler, \emph{On complex affine surfaces with $\mathbb{C}_+$-actions},
  Comment. Math. Helvetici \textbf{69} (1994), 5--27.

\bibitem{Giz71}
M.H. Gizatullin, \emph{Quasihomogeneous affine surfaces}, Math. USSR Izvestiya
  \textbf{5} (1971), 1057--1081.

\bibitem{EGAI}
A.~Grothendieck and J.~Dieudonné, \emph{{EGA I}. {L}e langage des sch{é}mas},
 Publ. {M}ath. {IHES}  vol.~4, 1960.

\bibitem{Ho72}
M.~Hochster, \emph{Nonuniqueness of coefficient rings in apolynomial ring},
  Proc. {A}mer. {M}ath. {S}oc. \textbf{34} (1972), 81--82.

\bibitem{Ii77}
S.~Iitaka, \emph{On logarithmic kodaira dimension of algebraic varieties},
  Complex analysis and algebraic geometry, {I}wanami {S}hoten, {T}okyo (1977),
  175--189.

\bibitem{IiF77}
S.~Iitaka and T.~Fujita, \emph{Cancellation theorem for algebraic varieties},
  J. {F}ac. {S}ci. {U}niv. {T}okyo \textbf{24} (1977), 123--127.

\bibitem{KML97}
S.~Kaliman and L.~Makar-Limanov, \emph{On the {R}ussel-{K}oras contractible
  threefolds}, J. Algebraic Geom. \textbf{6} (1997), 247--268.

\bibitem{ML96}
L.~Makar-Limanov, \emph{On the hypersurface $x+x^2y+z^2+t^3=0$ in
  $\mathbb{C}^4$ or a $\mathbb{C}^3$-like threefold which is not
  $\mathbb{C}^3$}, Israel J. Math. \textbf{96} (1996), 419--429.

\bibitem{ML01}
L.~Makar-Limanov, \emph{On the group of automorphisms of a surface
  $x^ny=p\left(z\right)$}, Israel J. Math. \textbf{121} (2001), 113--123.

\bibitem{MiyBook}
M.~Miyanishi, \emph{Open {A}lgebraic {S}urfaces}, vol.~12, CRM {M}onograph
  {S}eries, 2001.

\bibitem{Wil98}
J.~Wilkens, \emph{On the cancellation problem for surfaces}, {C}. {R}. {A}cad.
  {S}ci. {P}aris {Sé}r. {I} {M}ath. \textbf{326} (1998), no.~9, 1111--1116.

\end{thebibliography}

\end{document}